%% file: run.tex
\documentclass{conm-p-l}

\newtheorem{theorem}{Theorem}[section]
\newtheorem{lemma}[theorem]{Lemma}

\theoremstyle{definition}
\newtheorem{definition}[theorem]{Definition}

\theoremstyle{remark}

\numberwithin{equation}{section}

\newcommand{\sq}{\Upsilon}

\newcommand{\tu}{\tilde{u}}

\newcommand{\F}{{\mathcal F}}
\newcommand{\e}{\epsilon}
\newcommand{\vth}{\vartheta}
\newcommand{\U}{{\mathcal U}}
\newcommand{\N}{{\mathcal N}}

\newcommand{\A}{{\mathcal A}}

\newcommand{\R}{{\mathcal R}}
\newcommand{\BD}{B\"acklund-Darboux transformations}
\renewcommand{\k}{\kappa}
\newcommand{\ga}{\gamma}

\newcommand{\dl}{\delta}
\newcommand{\Dl}{\Delta}
\renewcommand{\th}{\theta}
\newcommand{\ra}{\rightarrow}
\newcommand{\al}{\alpha}
\newcommand{\be}{\beta}
\newcommand{\sg}{\sigma}
\newcommand{\Sg}{\Sigma}

\newcommand{\pa}{\partial}

\newcommand{\hq}{\hat{q}}

\newcommand{\hy}{\hat{y}}
\newcommand{\hz}{\hat{z}}

\newcommand{\la}{\lambda}
\newcommand{\bq}{\bar{q}}

\newcommand{\om}{\omega}
\newcommand{\Om}{\Omega}

\newcommand{\lag}{\langle}
\newcommand{\rag}{\rangle}
\newcommand{\tx}{\tilde{x}}

\renewcommand{\O}{{\mathcal O}}

\newcommand{\tz}{\tilde{z}}

\newcommand{\htau}{\hat{\tau}}
\newcommand{\hrho}{\hat{\rho}}
\newcommand{\hvth}{\hat{\vartheta}}
\newcommand{\tb}{\tilde{b}}

\begin{document}

\title{Melnikov Analysis for Singularly Perturbed DSII Equation}

\author{Yanguang (Charles)  Li}
\address{Department of Mathematics, University of Missouri, 
Columbia, MO 65211}
\curraddr{}
\email{cli@math.missouri.edu}
\thanks{}


\date{}



\begin{abstract}
Rigorous Melnikov analysis is accomplished for Davey-Stewartson II 
equation under singular perturbation. Unstable fiber theorem and 
center-stable manifold theorem are established. The fact that the 
unperturbed homoclinic orbit, obtained via a Darboux transformation, 
is a classical solution, leads to the conclusion that only local 
well-posedness is necessary for a Melnikov measurement.  
\end{abstract}

\maketitle








\input intro.tex

\input local.tex

\input global.tex

\input appendix.tex

\end{document}

%% file: intro.tex
\section{Introduction}

To build a Melnikov analysis for high dimensional nonlinear wave equations, 
we consider the Davey-Stewartson II equation (DSII) under a singular 
perturbation
\begin{equation}
\left \{ \begin{array}{l} iq_t=\sq q+ [2(|q|^2-\omega^2)+ 
u_y ]q +i\epsilon [\Delta q-\alpha q+\beta ]\ , \cr
\Delta u = -4\partial_y |q|^2 \ , \cr \end{array} \right.
\label{PDS}
\end{equation}
where $q$ is a complex-valued function of the three variables ($t,x,y$), 
$u$ is a real-valued function of the three variables ($t,x,y$), 
the external parameters $\omega$, $\alpha$, and $\beta$ are all 
positive constants, and $\epsilon>0$ is the perturbation parameter,
\[
\sq =\partial_{xx}-\partial_{yy}, \quad 
\Delta=\partial_{xx}+\partial_{yy}, \quad i=\sqrt{-1}.
\]
Periodic boundary condition is imposed,
\[
q(t,x+2\pi/\k_1,y)=q(t,x,y)=q(t,x,y+2\pi/\k_2),
\]
\[
u(t,x+2\pi/\k_1,y)=u(t,x,y)=u(t,x,y+2\pi/\k_2),
\]
where $\k_1$ and $\k_2$ are positive constants. Even constraint is also 
imposed,
\[
q(t,-x,y)=q(t,x,y)=q(t,x,-y),
\]
\[
u(t,-x,y)=u(t,x,y)=u(t,x,-y).
\]
Further constraints are placed upon $\om$, $\al$, $\be$, $\k_1$, and $\k_2$.
The first one $0<\al \om <\be$ is the condition for the existence of a 
saddle, and the second one is the condition for the existence of only 
two unstable modes,
\begin{equation}
\left \{ \begin{array}{l} \k_2 < \k_1 < 2\k_2 , \cr 
\k_1^2 < 4 \om^2 < \min \{ \k_1^2 +\k_2^2 , 4 \k_2^2 \} , 
\cr \end{array} \right. 
\label{cstr1}
\end{equation}
or
\begin{equation}
\left \{ \begin{array}{l} \k_1 < \k_2 < 2\k_1 , \cr 
\k_2^2 < 4 \om^2 < \min \{ \k_1^2 +\k_2^2 , 4 \k_1^2 \} . 
\cr \end{array} \right. 
\label{cstr2}
\end{equation}
Davey-Stewartson II equation can be regarded as a generalization of the 
1D cubic nonlinear Schr\"odinger equation (NLS) \cite{Li01b}. In fact, it 
is a nontrivial generalization in the sense that the spatial part of the 
Lax pair of the DSII is a system of two first order partial differential 
equations, for which there is no Floquet discriminant to describe the 
isospectral property, in contrast to the case for NLS. It turns out that 
Melnikov vectors can still be obtained through quadratic products of 
Bloch eigenfunctions, instead of the gradient of the Floquet discriminant 
as in the NLS case. 

At the moment, there is no global well-posedness for DSII in Sobolev 
spaces. In fact, DSII has finite-time blow-up solutions in 
$H^s(R^2)$, ($0<s<1$), \cite{Oza92}. Of course, DSII has local well-posedness 
in Sobolev spaces \cite{GS90} \cite{Sun94}. The Melnikov measurement 
is built upon an unperturbed homoclinic orbit of the unperturbed DSII. 
Explicit expression of such a homoclinic orbit can be obtained through 
Darboux transformation \cite{Li00}. The homoclinic orbit is a classical 
solution. This enables us to iterate the local well-posedness result in 
time, and complete a Melnikov measurement. Unstable fiber theorem and 
center-stable manifold theorem are of course needed, and established 
along the same line as in \cite{Li01b}. Novelties in regularity are 
introduced by the singular perturbation $\e \Dl$ which generates the 
semigroup $e^{\e t \Dl}$. 

The article is organized as follows: section 2 deals with local theory 
which includes unstable fiber theorem and center-stable manifold 
theorem, and we handle global theory in section 3 which includes integrable
theory and Melnikov analysis.

%% file: local.tex
\section{Local Theory}

One can view the perturbed DSII (\ref{PDS}) as an evolution equation in 
the $q$ variable. First, one can define the spatial mean as
\[
\lag q \rag = \frac{\k_1 \k_2}{4 \pi^2} \int_0^{2 \pi/\k_2}
\int_0^{2 \pi/\k_1} q \ dx dy .
\]
Then one may introduce the space $\dot{H}^s$ as 
\[
\dot{H}^s = \{ q \in H^s \ | \ \lag q \rag = 0 \}.
\]
The inverse Laplacian $\Dl^{-1} : \dot{H}^s \mapsto \dot{H}^{s+2}$ is an 
isomorphism. The perturbed DSII (\ref{PDS}) can be rewritten as
\begin{equation}
iq_t=\sq q+2[\Delta^{-1}\sq |q|^2+ \lag|q|^2 \rag -\omega^2]q
+i\epsilon (\Delta q-\alpha q+\beta ). 
\label{PDS2}
\end{equation}

\subsection{Change of Coordinates}

Denote by $\Pi$ the 2D subspace
\begin{equation}
\Pi = \{ q \ | \ \pa_x q = \pa_y q = 0 \} .
\label{Pi}
\end{equation}
Dynamics in $\Pi$ is the same as that given in \cite{Li01b}. Denote by 
$S_\om$ the circle
\begin{equation}
S_\om = \{ q \in \Pi \ | \ |q|=\om \}.
\label{circle}
\end{equation}
When $\al \om < \be$, there is a saddle $Q_\e$ near $S_\om$ in $\Pi$, 
which is located at $q= I e^{i\th}$ where
\begin{equation}\label{Qec}\begin{cases}I=\omega^2-\epsilon 
\frac{1}{2\omega}\sqrt{\beta^2-\alpha^2\omega^2}+\cdots ,\\
\cos \theta  =\frac{\alpha \sqrt{I}}{\beta},&\theta \in \left( 
0,\frac{\pi}{2}\right).\end{cases}\end{equation}
Its eigenvalues are
\begin{equation}\label{Qee} \mu_{1,2}=\pm 
\sqrt{\epsilon}\sqrt{4\sqrt{I}\beta\sin \theta -\epsilon \left( 
\frac{\beta \sin
\theta}{\sqrt{I}}\right)^2}-\epsilon \alpha,\end{equation}
where $I$ and $\theta $ are given in \eqref{Qec}.
In the entire phase space, $Q_\e$ is still a saddle. Local theory 
will be built in a tubular neighborhood of $S_\om$. Let
\[
q(t,x,y)=[\rho(t) +f(t,x,y)]e^{i\theta(t)},\quad \lag f \rag = 0.
\]
Let 
\[
I = \lag |q|^2 \rag = \rho^2 + \lag |f|^2 \rag , \quad  J=I-\om^2.
\]
In terms of the new variables ($J, \th, f$), Equation (\ref{PDS2}) can 
be rewritten as 
\begin{eqnarray}
\dot{J} &=& \epsilon \bigg [ -2\alpha(J+\omega^2)+2\beta\sqrt{J+\omega^2}
\cos \theta \bigg ] +\epsilon \R_2^J, \label{cc1} \\
\dot{\th} &=& -2J - \epsilon \beta \frac {\sin \theta}{\sqrt{J+\omega^2}}
+\R_2^\th, \label{cc2} \\
f_t &=& L_\epsilon f+V_\epsilon f-i \N_2 -i \N_3, \label{cc3}
\end{eqnarray} 
where 
\begin{eqnarray*}
L_\epsilon f &=& -i\sq f+\epsilon (\Delta-\alpha )f-2i\omega^2
\Delta^{-1}\sq (f+\bar f), \\
V_\epsilon f &=& -i2J\Delta^{-1}\sq (f+\bar f)
+i\epsilon \beta \frac {\sin \theta}{\sqrt{J+\omega^2}}f, \\
\R^J_2 &=& -2 \lag \nabla f  \cdot \nabla \bar{f} \rag + 2\beta \cos \theta
\bigg [ \sqrt{J+\omega^2-\lag f \rag^2}-\sqrt{J+\omega^2} \bigg ], \\
\R^\th_2 &=& - \lag (f +\bar{f} ) \Dl^{-1} \sq (f+\bar{f} )\rag 
-\rho^{-1} \lag (f+\bar{f})\Delta^{-1}\sq |f|^2 \rag  \\
& & -\epsilon \beta \sin \theta \bigg [ \frac {1}{\sqrt{J+\omega^2-
\lag |f|^2\rag }} -\frac {1}{\sqrt{J+\omega^2}}\bigg ], \\
\N_2 &=& 2\rho \bigg [ \Delta^{-1}\sq |f|^2
+f\Delta^{-1} \sq (f+\bar f) 
- \lag f\Delta^{-1}\sq (f+\bar f)\rag \bigg ],\\
\N_3 &=& -\lag (f+\bar f)\Delta^{-1}\sq(f+\bar f) \rag f 
- 2 \lag |f|^2 \rag \Delta^{-1}\sq(f+\bar f) \\
& & +2 \bigg [ f\Delta^{-1}\sq|f|^2 - \lag 
f\Delta^{-1}\sq|f|^2 \rag \bigg ] 
- \rho^{-1} \lag (f+\bar f)\Delta^{-1}\sq|f|^2 \rag f \\
& & - \epsilon \beta \sin \theta \bigg [ \frac {1}
{\sqrt{J+\omega^2-\lag |f|^2 \rag }}
-\frac {1}{\sqrt{J+\omega^2}} \bigg ]f.
\end{eqnarray*}
Since $H^s$ ($s \geq 2$) is a Banach algebra \cite{Ada75}, we have
\begin{eqnarray*}
& & |\R^J_2| \sim  \mathcal{O}(\| f\|^2_s),\quad |\mathcal{R}^\theta_2|
\sim \mathcal{O}(\|f\|^2_s), \\
& & \|\mathcal{N}_2\|_s\sim \mathcal{O} (\| f\|^2_s),
\quad \| \mathcal{N}_3\|_s\sim \mathcal{O}(\| f\|^3_s),
\quad (s\geq 2). 
\end{eqnarray*}

\subsection{Unstable Fibers}

On $\Pi$ (\ref{Pi}), the saddle $Q_\e$ has an unstable and a stable curves 
which lie in an annular neighborhood of $S_\om$ in $\Pi$. The width of this 
annular neighborhood is of order $\O(\sqrt{\e})$. 
\begin{definition}
For any $\dl > 0$, we define the annular neighborhood of the circle
$S_\om$ (\ref{circle}) in $\Pi$ (\ref{Pi}) as 
\begin{equation}
\A(\dl) = \{ (J, \th)\ | \ |J| < \dl \}.
\label{ann}
\end{equation}
\end{definition}
Unstable fibers with base points in $\A(\dl)$ for some $\dl > 0$ persist,
even under the singular perturbation.

The spectrum of $L_\e$ consists of only point spectrum. The eigenvalues of 
$L_\e$ are:
\begin{equation}
\mu_\xi^{\pm} = -\e (\al +|\xi|^2)\pm |\xi|^{-1}|\xi_1^2 -\xi_2^2| 
\sqrt{4\om^2 - |\xi|^2}, 
\label{eval}
\end{equation}
where $\xi =(\xi_1,\xi_2)$, $\xi_j = k_j \k_j$, $k_j =0,1,2,\cdots$,
($j=1,2$), $k_1+k_2 > 0$, $|\xi|^2 =\xi_1^2+\xi_2^2$, and $\k_1$, $\k_2$, 
and $\om$ satisfy the constraint (\ref{cstr1}) or (\ref{cstr2}).

Denote $\mu_{(\k_1,0)}^\pm$ by $\mu_x^\pm$ and $\mu_{(0,\k_2)}^\pm$ by 
$\mu_y^\pm$. The eigenfunctions corresponding to $\mu_x^\pm$ and 
$\mu_y^\pm$ are
\[
u_x^\pm=e^{\pm i \vth_x} \cos \k_1 x\ , \quad e^{\pm i \vth_x} =
\frac{\k_1 \mp i \sqrt{4\om^2 -\k_1^2}}{2\om}\ , 
\]
\[
u_y^\pm=e^{\pm i \vth_y} \cos \k_2 y\ , \quad e^{\pm i \vth_y} =
\frac{\k_2 \pm i \sqrt{4\om^2 -\k_2^2}}{2\om}\ . 
\]
Notice also that the singular perturbation $-\e |\xi|^2$ breaks the gap 
between the center spectrum and the stable spectrum. Nevertheless, the 
gap between the unstable spectrum and the center spectrum survives. This 
leads to the following unstable fiber theorem.
\begin{theorem}[Unstable Fiber Theorem]
For any $s \geq 2$, there exists a $\dl >0$ such that for any 
$p \in \A(\dl)$, there is an 
unstable fiber $\F^u_p$ which is a 2D surface. $\F^u_p$ has
the following properties:
\begin{enumerate}
\item $\mathcal{F}^u_p$ is a $C^1$ smooth surface in 
$\|  \ \|_s$ norm. 
\item $\mathcal{F}^u_p$ is also $C^1$ smooth in $\epsilon$, $\alpha$, 
$\beta$, $\omega$, and $p$ in $\|\ \|_s$ norm, 
$\epsilon \in [0,\epsilon_0)$ for some $\epsilon_0>0$ depending on $s$.
\item $p\in \mathcal{F}^u_p$, $\mathcal{F}^u_p$ is tangent to $\ \mbox{span}\  
\{ u^+_x,u^+_y\}$ at $p$ when $\epsilon=0$.
\item $\mathcal{F}^u_p$ has the exponential decay property: Let $S^t$ 
be the evolution operator of \eqref{cc1}-\eqref{cc3}, $\forall
p_1\in \mathcal{F}^u_p$,
\[
\| S^tp_1-S^tp\|_s \leq 
Ce^{\frac{1}{3}\mu^+t} \| p_1-p\|_s ,\quad \forall t\leq 0,
\]
where $\mu^+=\min \{ \mu^+_x,\mu^+_y\}$.
\item $\{ \mathcal{F}^u_p\}_{p\in \mathcal{A}(\dl)}$ forms an invariant 
family of unstable fibers,
\[
S^t\mathcal{F}^u_p\subset \mathcal{F}^u_{S^tp}\ ,\quad 
\forall t\in [-T,0],
\]
and $\forall T>0$ ($T$ can be $+\infty$), such that $S^\tau p\in 
\mathcal{A}(\dl)$, $\forall \tau \in [-T,0]$.
\end{enumerate}
\label{UFT} 
\end{theorem}
The proof of this theorem follows from the same arguments as in \cite{Li01b}.
Notice, in particular, that $\F^u_p \subset H^s$ for any $s \geq 2$. It 
is this fact that leads to the $C^1$ smoothness of $\F^u_p$ in $\e$. 
Denote by $W^u(Q_\e)$ the unstable manifold of the saddle $Q_\e$ (\ref{Qec}),
which is 3-dimensional. Denote by $W^u_\Pi(Q_\e)$ the unstable curve of 
$Q_\e$ in $\Pi$ (\ref{Pi}). $W^u_\Pi(Q_\e)=\Pi \cap W^u(Q_\e)$, and 
$W^u_\Pi(Q_\e) \subset \A(\dl)$. $W^u(Q_\e)$ has the fiber representation 
\begin{equation}
W^u(Q_\e)=\bigcup_{p \in W^u_\Pi(Q_\e)} \F^u_p \ .
\label{uqm}
\end{equation}
Thus $W^u(Q_\e) \subset H^s$ for any $s \geq 2$.

\subsection{Center-Stable Manifold}

Also due to the fact that the gap between unstable spectrum and center 
spectrum survives under the singular perturbation (\ref{eval}), a 
center-stable manifold persists.
\begin{theorem}[Center-Stable Manifold Theorem]
There exists a $C^1$ smooth codimension 2 
locally invariant center-stable manifold $W^{cs}_n$ in $H^n$ for any
$n\geq 2$.
\begin{enumerate}
\item At points in the subset $W^{cs}_{n+4}$ of 
$W^{cs}_n$, $W^{cs}_n$ is $C^1$ smooth in $\epsilon$, in $H^n$ norm, 
for $\epsilon \in
[0,\epsilon_0)$ and some $\epsilon_0 >0$. 
\item $W^{cs}_n$ is $C^1$ smooth in $(\alpha ,\beta ,\omega)$.
\item The annular neighborhood $\mathcal{A}(\dl)$ in Theorem~\ref{UFT} is 
included in $W^{cs}_n$.
\end{enumerate}
\label{CSM}
\end{theorem}
The proof of this theorem follows from the same arguments as in \cite{Li01b}.

Regularity of $W^{cs}_n$ in $\e$ is crucial in Melnikov analysis. Melnikov 
integrals are the leading order terms in $\e$ of the signed distances between 
$W^u(Q_\e)$ (\ref{uqm}) and $W^{cs}_n$. The signed distances are set up along
an unperturbed homoclinic orbit, and the regularity of $W^{cs}_n$ in $\e$ 
at $\e = 0$ determines the order of the signed distances in $\e$. Due to 
the singular perturbation, $W^{cs}_n$ is not $C^1$ in $\e$ at every point
rather at points in the subset $W^{cs}_{n+4}$. Here one may be able to 
replace $W^{cs}_{n+4}$ by $W^{cs}_{n+2}$. But we are not interested in 
sharper results, and the current result is sufficient for our purpose.

\subsection{Local Well-Posedness}

Following a much easier argument than that in e.g. \cite{Kat72} \cite{Kat75},
one can prove the following local well-posedness theorem.
\begin{theorem}
For any $q_0 \in H^n$ ($n \geq 2$), there exists $\tau =\tau (\|q_0\|_n) > 0$, 
such that the perturbed DSII (\ref{PDS2}) has a unique solution 
$q(t) = S^t(q_0; \e, \al, \be, \om) \in C^0([0,\tau ],H^n)$, $q(0)=q_0$, 
where $S^t$ denotes the evolution operator. $S^t(\cdot ; \e, \al, \be, \om):
H^n \mapsto H^n$ is $C^1$ in $q_0$ and ($\al,\be ,\om$). 
$S^t(\cdot ; \e, \al, \be, \om): H^{n+4} \mapsto H^n$ is $C^1$ in
$t$ and $\e$, $\e \in [0, \e_0)$, $\e_0 >0$.
\label{locwp}
\end{theorem}
Here ``$C^1$ in $q_0$ and ($\al,\be ,\om$)'' can be replaced by 
``$C^\infty$ in $q_0$ and ($\al,\be ,\om$)''. $H^{n+4}$ can be replaced by 
$H^{n+2}$. But we are not interested in sharper results.

%% file: global.tex
\section{Global Theory}

Global theory is referred to a theory global in phase space, which includes 
integrable theory and Melnikov analysis. Integrable theory provides two 
ingredients for a Melnikov analysis: (1). An explicit expression of the 
unperturbed homoclinic orbit, (2). Melnikov vectors with explicit expressions.

\subsection{Integrable Theory}

Calculations in this subsection are essentially the same with those in 
\cite{Li00}. The minor differences are introduced by the spatial periods
$2\pi/\k_1$ and $2\pi/\k_2$ in contrast to $2\pi$ and $2\pi$ in \cite{Li00}.
Proofs of theorems and lemmas can be found in \cite{Li00}.

The DSII [$\e =0$ in (\ref{PDS})] is an integrable system with the Lax 
pair
\begin{eqnarray}
         L \psi &=& \lambda \psi\,, \label{LP1} \\
       \partial_t \psi &=& A \psi\,,\label{LP2}
\end{eqnarray}
where $\psi = \left( \psi_1, \psi_2\right)^T$, and 
\[
  L = \left(
\begin{array}{lr}
D^{-} & q\\ \\
\bq & D^{+}
\end{array}
\right)\,,
\]
\[
A = i \left[
2 \left(
\begin{array}{cc}
- \partial^2_x & q \partial_x\\
\bq \partial_x & \partial^2_x
\end{array}
\right) \, + \,
\left(
\begin{array}{cc}
r_1 & (D^+ q)\\
-(D^{-} \bq ) & r_2
\end{array}
\right)
\right]\, ,
\]
where
\begin{equation}
  D^+ = \alpha \partial_y + \partial_x\,, \qquad D^{-} = \alpha
  \partial_y - \partial_x\, , \qquad \al^2 = -1\ ,
\label{DD}
\end{equation}
$r_1$ and $r_2$ have the expressions,
\begin{equation} 
  r_1 = \frac{1}{2} [-2(|q|^2-\om^2)-u_y + i \tu ] \ , \quad r_2= 
\frac{1}{2} [2(|q|^2-\om^2)+u_y + i \tu ] \, , 
\label{res1}
\end{equation}
where $\tu$ is also a real-valued function satisfying
\[
\Dl \tu = 4i\al \pa_x \pa_y |q|^2\ .
\]
Notice that DSII is invariant under the transformation $\sg$:
\begin{equation}
  \sigma (q,\bq, r_1, r_2; \al) = (q,\bq, -r_2, -r_1; -\al)\,.
\label{IT}
\end{equation}
Applying the transformation $\sigma$ (\ref{IT}) to the Lax pair
(\ref{LP1}, \ref{LP2}), we have a congruent Lax pair for which
the compatibility condition gives the same DSII.
The congruent Lax pair is given as:
\begin{eqnarray}
 \hat{L} \hat{\psi} &=& \lambda \hat{\psi}
                               \,,\label{CLP1} \\
 \partial_t \hat{\psi} &=& \hat{A}
                                   \hat{\psi}
                                 \,,\label{CLP2}
\end{eqnarray}
where $\hat{\psi} = (\hat{\psi}_1, \hat{\psi}_2)$, and
\[
  \hat{L} =
  \left(
    \begin{array}{cc}
- D^+ & q\\ \\
\bq & -D^-
    \end{array}
  \right)\,,
\]
\[
\hat{A} = i \left[
  2 \left(
    \begin{array}{cc}
- \partial^2_x & q \partial_x\\
\bq \partial_x & \partial^2_x
    \end{array}
  \right) + 
  \left(
    \begin{array}{cc}
-r_2 & -(D^- q)\\
(D^+\bq ) & -r_1
    \end{array}
  \right)
\right]\,.
\]

The B\"acklund-Darboux transformation can be formulated as follows. 
Let $(q,u)$ be a solution to the DSII, and let $\lambda_*$ be any 
value of $\lambda$. Let
$\psi = (\psi_1, \psi_2)^T$ be a solution to the Lax
pair (\ref{LP1}, \ref{LP2}) at $(q, \bar q, r_1, r_2;
\lambda_*)$. Define the matrix operator:
\begin{displaymath}
  \Gamma = 
\left[
  \begin{array}{cc}
             \wedge + a & b\\
             c & \wedge + d
  \end{array}
\right]\,,
\end{displaymath}
where $\wedge = \alpha \partial_y - \lambda$, and $a$, $b$, $c$,
$d$ are functions defined as:
\begin{eqnarray*}
  a &=& \frac{1}{\Delta} \left[ \psi_2 \wedge_2 \bar{\psi}_2 +
                  \bar{\psi}_1 \wedge_1 \psi_1 \right]\,,\\[2ex]
  b &=& \frac{1}{\Delta} \left[ \bar{\psi}_2 \wedge_1 \psi_1 -
                   \psi_1 \wedge_2 \bar{\psi}_2 \right]\,,\\[2ex]
  c &=& \frac{1}{\Delta} \left[ \bar{\psi}_1 \wedge_1 \psi_2
                     - \psi_2 \wedge_2 \bar{\psi}_1 \right]\,,\\[2ex]
  d &=&  \frac{1}{\Delta} \left[ \bar{\psi}_2 \wedge_1 \psi_2 +
                      \psi_1 \wedge_2 \bar{\psi}_1 \right]\,,
\end{eqnarray*}
in which $\wedge_1 = \alpha \partial_y - \lambda_*$, $\wedge_2 =
\alpha \partial_y + \overline{\lambda_*}$, and
\begin{displaymath}
  \Delta = - \left[ | \psi_1 |^2 + |\psi_2|^2 \right]\,.
\end{displaymath}
Define a transformation as follows:
\begin{displaymath}
  \left\{
    \begin{array}{ccc}
(q,r_1,r_2) &\rightarrow& (Q,R_1,R_2)\,, \\
\phi &\rightarrow& \Phi\,;
    \end{array}
\right.
\end{displaymath}
\begin{eqnarray}
                    Q   &=& q - 2b\,,\nonumber \\[2ex]
                    R_1 &=& r_1 + 2(D^+a)\,, \label{DSBT}\\[2ex]
                    R_2 &=& r_2 - 2 (D^- d)\,,\nonumber\\[2ex]
                   \Phi &=& \Gamma \phi\,;\nonumber
\end{eqnarray}
where $\phi$ is any solution to the Lax pair (\ref{LP1},
\ref{LP2}) at $(q, \bar{q}, r_1, r_2; \lambda)$, $D^+$
and $D^-$ are defined in (\ref{DD}), we have the following theorem 
\cite{Li00}.
\begin{theorem}
The transformation (\ref{DSBT}) is a B\"acklund-Darboux
transformation. That is, the function $Q$ defined
through the transformation (\ref{DSBT}) is also a solution to the 
DSII. The function $\Phi$ defined through the transformation
(\ref{DSBT}) solves the Lax pair (\ref{LP1}, \ref{LP2}) at $(Q,
\bar{Q}, R_1, R_2; \lambda)$.
\label{DSTH}
\end{theorem}

Consider the spatially independent solution,
\begin{equation}
q_c = \eta \exp \{ -2i [ \eta^2 - \om^2 ] t + i \ga \} \ ,
\label{us}
\end{equation}
where $\eta$ satisfies the constraint (\ref{cstr1}) and (\ref{cstr2})
with $\om$ replaced by $\eta$.
The dispersion relation for the linearized DSII at $q_c$ is 
\[
\Om = \pm \frac{|\xi_1^2 - \xi_2^2|}{\sqrt{\xi_1^2 +\xi_2^2}}
\sqrt{4 \eta^2 - (\xi_1^2 +\xi_2^2)}\ , \ \ \mbox{for} \ 
\dl q \sim q_c \exp \{ i (\xi_1 x +\xi_2 y) +\Om t \} \ ,
\]
where $\xi_1 = k_1 \k_1$, $\xi_2 = k_2 \k_2$, and $k_1$ and $k_2$ 
are integers. There are only 
two unstable modes ($\k_1, 0$) and ($0, \k_2$) under even constraint.

The Bloch eigenfunction of the Lax pair (\ref{LP1}) and (\ref{LP2})
is given as,
\begin{equation}
\tilde{\psi} = c(t) \left[
    \begin{array}[]{c}
-q_c \\ \chi
    \end{array} \right]
\exp \left\{ i(\xi_1 x + \xi_2y) \right\} \, ,
\label{slLax}
\end{equation}
where
\begin{eqnarray*}
& & c(t) = c_0 \exp \left\{ \left[ 2\xi_1(i \alpha \xi_2 - \lambda )
            + ir_2 \right] t \right\} \, , \\
& & r_2 - r_1 = 2 ( \left| q_c \right|^2 - \omega^2 ) \, , \\
& & \chi = (i \alpha \xi_2- \lambda )-i\xi_1 \, , \\
& & (i \alpha \xi_2 - \lambda)^2 + \xi^2_1 = \eta^2 \, .
\end{eqnarray*}
For the iteration of the B\"acklund-Darboux transformations, one 
needs two sets of eigenfunctions. First, we choose
$\xi_1 = \pm \frac{1}{2} \k_1$, $\xi_2=0$, $\lambda_0 = \sqrt{\eta^2 -
\frac{1}{4}\k_1^2}$ (for a fixed branch),
\begin{eqnarray}
  \psi^{\pm} = c^{\pm} \left[
    \begin{array}{c}
      -q_c \\ \\ \chi^{\pm}
    \end{array} \right]
  \exp \left\{ \pm i \frac{1}{2}\k_1x \right\} \, , 
\label{efunc1}
\end{eqnarray}
where
\begin{eqnarray*}
& & c^{\pm} = c^{\pm}_0 \exp \left\{ \left[ \mp \k_1 \lambda_0 + ir_2
    \right] t \right\} \, , \\
& & \chi^{\pm} = - \lambda_0 \mp i \frac{1}{2} \k_1 =
  \eta e^{\mp i (\frac{\pi}{2} +\vth_1)} \, , \ \mbox{i.e.} \ 
\eta e^{\pm i \vth_1}= \frac{1}{2} \k_1 \pm i\lambda_0 \ .
\end{eqnarray*}
We apply the B\"acklund-Darboux transformations with $\psi = 
\psi^+ + \psi^-$, which generates
the unstable foliation associated with the $(\k_1,0)$ 
linearly unstable mode. Then, we choose $\xi_2 = \pm
\frac{1}{2}\k_2$, $\lambda =0$, $\xi^0_1 = \sqrt{\eta^2-\frac{1}{4}\k_2^2}$
(for a fixed branch),
\begin{equation}
  \phi_{\pm} = c_{\pm} \left[
    \begin{array}{c}
      -q_c \\ \\ \chi_{\pm}
    \end{array} \right] 
  \exp \left\{ i (\xi^0_1 x \pm \frac{1}{2}\k_2 y) \right\} \, , 
\label{efunc2}
\end{equation}
where
\begin{eqnarray*}
& & c_{\pm} = c^0_{\pm} \exp \left\{ \left[ \pm i \alpha \k_2\xi^0_1
      + ir_2 \right] t \right\} \, , \\
& & \chi_{\pm} = \pm i \alpha \frac{1}{2}\k_2 - i\xi^0_1 
=\pm \eta e^{\mp i \vth_2}\, , \ \mbox{i.e.}\ \eta e^{\pm i \vth_2} =
i \alpha \frac{1}{2}\k_2 \pm i\xi^0_1\ .
\end{eqnarray*}
We start from these eigenfunctions $\phi_{\pm}$ to generate
$\Gamma \phi_{\pm}$ through \BD, and then iterate the \BD ~with
$\Gamma \phi_+ + \Gamma \phi_-$ to generate the unstable
foliation associated with all the linearly unstable modes $(
\k_1,0)$ and $(0, \k_2)$. It turns out that the following 
representations are convenient,
\begin{eqnarray}
\psi^\pm &=& \sqrt{c_0^+c_0^-}e^{ir_2 t}\left ( \begin{array}{c} 
v_1^\pm \cr v_2^\pm \cr  \end{array} \right ) \ ,
\label{rwf1} \\
\phi_\pm &=& \sqrt{c^0_+c^0_-}e^{i \xi_1^0 x + ir_2 t}\left ( \begin{array}{c} 
w_1^\pm \cr w_2^\pm \cr \end{array} \right ) \ ,
\label{rwf2} 
\end{eqnarray}
where
\[
v_1^\pm = -q_c e^{\mp \frac{\tau}{2} \pm i \tx} \ , \ \ 
v_2^\pm = \eta e^{\mp \frac{\tau}{2} \pm i \tz} \ ,
\]
\[
w_1^\pm = -q_c e^{\pm \frac{\htau}{2} \pm i \hy} \ , \ \ 
w_2^\pm = \pm \eta e^{\pm \frac{\htau}{2} \pm i \hz}\ ,
\]
and
\[
c_0^+/c_0^- = e^{\rho + i \vth }\ , \ \ 
\tau = 2\k_1 \la_0 t - \rho \ , \ \ 
\tx = \frac{1}{2} \k_1 x + \frac{\vth}{2}\ ,  \ \ 
\tz = \tx - \frac{\pi}{2} - \vth_1 \ ,
\]
\[
c^0_+/c^0_- = e^{\hrho + i \hvth }\ , \ \ 
\htau = 2i\al \k_2 \xi_1^0 t + \hrho \ , \ \ 
\hy = \frac{1}{2} \k_2 y + \frac{\hvth}{2}\ , \ \
\hz = \hy - \vth_2 \ .
\]
The following representations are also very useful,
\begin{eqnarray}
\psi &=& \psi^+ + \psi^- = 2 \sqrt{c_0^+c_0^-}e^{ir_2 t}
\left ( \begin{array}{c} 
v_1 \cr v_2 \cr  \end{array} \right ) \ ,
\label{rwf3} \\
\phi &=& \phi^+ + \phi^- = 2 \sqrt{c^0_+c^0_-}e^{i \xi_1^0 x + 
ir_2 t}\left ( \begin{array}{c} 
w_1 \cr w_2 \cr \end{array} \right ) \ ,
\label{rwf4} 
\end{eqnarray}
where
\[
v_1 = -q_c [ \cosh \frac{\tau}{2} \cos \tx - i \sinh \frac{\tau}{2} \sin \tx ]
\ , \ \ v_2 = \eta [ \cosh \frac{\tau}{2} \cos \tz - i \sinh \frac{\tau}{2} 
\sin \tz ] \ , 
\]
\[
w_1 = -q_c [ \cosh \frac{\htau}{2} \cos \hy + i \sinh \frac{\htau}{2} \sin 
\hy ]
\ , \ \ w_2 = \eta [ \sinh \frac{\htau}{2} \cos \hz + i \cosh \frac{\htau}{2} 
\sin \hz ] \ . 
\]
Applying the \BD ~(\ref{DSBT}) with $\psi$ given in (\ref{rwf3}), we
have the representations,
\begin{eqnarray}
a &=& -\lambda_0 \ \mbox{sech}\ \tau \sin (\tx + \tz) \sin (\tx -\tz)
\nonumber \\
& & \quad \times \bigg [ 1 + \ \mbox{sech}\ \tau \cos (\tx + \tz) 
\cos (\tx -\tz) \bigg ]^{-1}\, , \label{aexp} \\
b &=& -q_c \tb = - \frac{\lambda_0 q_c}{\eta} \bigg [ \cos (\tx - \tz) -
i \tanh \tau \sin (\tx - \tz)
\nonumber \\
& & \quad +\ \mbox{sech}\ \tau \cos (\tx + \tz)\bigg ]
\bigg [ 1+ \ \mbox{sech}\ \tau  \cos (\tx + \tz) 
\cos (\tx - \tz)\bigg ]^{-1} \, , \label{bexp}\\
& & c= \overline{b} \, , \ \ \ \ d= - \overline{a} =-a \, . \label{cdexp}
\end{eqnarray}
The evenness of $b$ in $x$ is enforced by the requirement that 
$\vth - \vth_1 = \pm \frac{\pi}{2}$, and
\begin{eqnarray}
a^{\pm} &=& \mp \lambda_0 \ \mbox{sech} \ \tau \cos \vth_1 \sin (\k_1 x)
\nonumber \\
& & \quad \times \bigg [ 1 \mp \ \mbox{sech}\ \tau \sin \vth_1 \cos (\k_1 x)
\bigg ]^{-1}\, , \label{eaexp} \\
b^{\pm} &=& -q_c \tb^{\pm} = - \frac{\lambda_0 q_c}{\eta} \bigg [ 
-\sin \vth_1 -
i \tanh \tau \cos \vth_1 
\nonumber \\
& & \quad \pm \ 
\mbox{sech}\ \tau \cos (\k_1 x) \bigg ] 
\bigg [ 1 \mp \ \mbox{sech}\ \tau  \sin \vth_1 \cos (\k_1 x)
\bigg ]^{-1} \, , \label{ebexp}\\
& & c= \overline{b} \, , \ \ \ \ d= - \overline{a} =-a \, . \label{ecdexp}
\end{eqnarray}
Notice also that $a^{\pm}$ is an odd function in $x$.  Under 
the above \BD, the
eigenfunctions $\phi_{\pm}$ (\ref{efunc2}) and $\phi$ are transformed into
\begin{equation}
\varphi^{\pm} = \Gamma \phi_{\pm} \, , \quad 
\varphi = \Gamma \phi = \Gamma \phi_+ + \Gamma \phi_-\ ,
\label{tphi}
\end{equation}
where
\begin{eqnarray*}
  \Gamma = \left[
    \begin{array}{cc}
\Lambda + a & b \\ \\ 
\overline{b} & \Lambda -a
    \end{array} \right] \, ,
\end{eqnarray*}
and $\Lambda = \alpha \partial_y - \lambda$ with $\lambda$
evaluated at $0$. Then
\[
\varphi^{\pm} =\sqrt{c_+^0 c_-^0} e^{i \xi_1^0 x + i r_2 t}\left[ 
\begin{array}{c} -q_c W_1^\pm \cr \cr \eta W_2^\pm \cr \end{array} 
\right ]\ ,
\]
where
\begin{eqnarray*}
W_1^\pm &=& [\pm i \frac{1}{2} \al \k_2 +a \pm \eta \tb e^{\mp i \vth_2}]
e^{\pm \frac{\htau}{2} \pm i \hy }\ ,  \\
W_2^\pm &=& \pm e^{\mp i \vth_2}[\pm i \frac{1}{2} \al \k_2 -a \pm \eta 
\bar{\tb} e^{\pm i \vth_2}]e^{\pm \frac{\htau}{2} \pm i \hy }\ ;
\end{eqnarray*}
\[
\varphi =  2 \sqrt{c_+^0 c_-^0} e^{i \xi_1^0 x + i r_2 t}\left[
    \begin{array}{c}
      -q_c W_1 \\ \\ \eta W_2
    \end{array} \right] \ , 
\]
where
\begin{eqnarray*}
W_1 &=& \cosh \frac{\htau}{2} [ a \cos \hy -\frac{1}{2} \al \k_2\sin \hy
+i \eta \tb \sin \hz ] \\
& & +\sinh \frac{\htau}{2} [ \frac{1}{2} i\al \k_2\cos \hy +ia \sin \hy
+\eta \tb \cos \hz ] \ , \\
W_2 &=& \cosh \frac{\htau}{2} [ -ia \sin \hz +\frac{1}{2} i\al \k_2\cos \hz
+\eta \bar{\tb} \cos \hy ] \\
& & +\sinh \frac{\htau}{2} [ -\frac{1}{2} \al \k_2\sin \hz -a \cos \hz
+i\eta \bar{\tb} \sin \hy ] \ .
\end{eqnarray*}
We generate the coefficients in the \BD ~
(\ref{DSBT}) with $\varphi$ (the iteration of the \BD),
\begin{eqnarray}
a^{(I)} &=& - \bigg [ W_2 (\al \pa_y \overline{W_2}) + 
\overline{W_1} (\al \pa_y W_1) \bigg ]\bigg [ |W_1|^2 
+|W_2|^2 \bigg ]^{-1}\ ,  \label{ria}\\
b^{(I)} &=& \frac{q_c}{\eta}\bigg [ \overline{W_2} (\al \pa_y W_1) - 
W_1 (\al \pa_y \overline{W_2}) \bigg ]\bigg [ |W_1|^2 
+|W_2|^2 \bigg ]^{-1}\ , \label{rib} \\
& & c^{(I)} = \overline{b^{(I)}} \, , \ \ \ \ \, d^{(I)} =-\overline{a^{(I)}} 
\, , \label{ricd}
\end{eqnarray}
where
\begin{eqnarray*}
& & W_2 (\al \pa_y \overline{W_2}) + \overline{W_1} (\al \pa_y W_1) \\
& & = \frac{1}{2} \al \k_2 \bigg \{ \cosh \htau \bigg [ -\al \k_2 a 
+ i a \eta (\tb + \overline{\tb}) \cos \vth_2 \bigg ]  \\
& & +\bigg [ \frac{1}{4}\k_2^2 - a^2 -\eta^2 |\tb|^2 \bigg ] 
\cos (\hy +\hz ) \sin \vth_2 + \sinh \htau \bigg [ a \eta  
(\tb - \overline{\tb}) \sin \vth_2 \bigg ]  \bigg \} \ , \\ \\
& & |W_1|^2 +|W_2|^2 \\
& & = \cosh \htau \bigg [ a^2 + \frac{1}{4}\k_2^2 + \eta^2 |\tb|^2  
+i\al \k_2 \eta \frac{1}{2} (\tb + \overline{\tb})\cos \vth_2 \bigg ] \\
& & +\bigg [ \frac{1}{4}\k_2^2 - a^2 -\eta^2 |\tb|^2 \bigg ]  
\sin (\hy +\hz ) \sin \vth_2 + \sinh \htau \bigg [ \al \k_2 \eta  
\frac{1}{2} (\tb - \overline{\tb}) \sin \vth_2 \bigg ]\ , \\ \\
& & \overline{W_2} (\al \pa_y W_1) - W_1 (\al \pa_y \overline{W_2}) \\ 
& & =\frac{1}{2} \al \k_2 \bigg \{ \cosh \htau \bigg [ -\al \k_2 \eta \tb  
+ i (-a^2 + \frac{1}{4}\k_2^2 + \eta^2 \tb^2 ) \cos \vth_2 \bigg ]  \\
& & + \sinh \htau \bigg [a^2 - \frac{1}{4}\k_2^2 + \eta^2 \tb^2 \bigg ]
\sin \vth_2 \bigg \} \ . 
\end{eqnarray*}
The new solution to the DSII is given by
\begin{equation}
   Q= q_c -2b -2b^{(I)} \, .
\label{newsl} 
\end{equation}
The evenness of $b^{(I)}$ in $y$ is enforced by the requirement
that $\hvth - \vth_2 = \pm \frac{\pi}{2}$.  In fact, we have
\begin{lemma}
Choosing the B\"acklund parameters $\vth$ and $\hvth$ as follows: 
$\vth = \vth_1  \pm \frac{\pi}{2}$, and
$\hvth = \vth_2  \pm \frac{\pi}{2}$,
\begin{equation}
  b(-x) = b(x) \, , \ \ \  b^{(I)}(-x,y)=b^{(I)}(x,y)=b^{(I)}(x,-y) \, , 
\end{equation}
and $Q = q_c -2b-2b^{(I)}$ is even in both $x$ and $y$.
\label{evenla}
\end{lemma}
The asymptotic behavior of $Q$ can be computed directly.  In
fact, we have the asymptotic phase shift lemma.
\begin{lemma}[Asymptotic Phase Shift Lemma]
For $\la_0 > 0$, $\xi_1^0 > 0$, and $\al = -i$; as $t \ra \pm \infty$, 
\begin{equation}
    Q = q_c -2b-2b^{(I)} \ra q_c e^{i\pi } e^{\mp i 2 (\vth_1 -\vth_2 )}\ .
\label{ayp}  
\end{equation}
In comparison, the asymptotic phase shift of the
first application of the \BD ~ is given by
\[
q_c - 2b \ra  q_c e^{\mp i 2 \vth_1}\ .
\]
\end{lemma}

Next we generate the Melnikov vectors.  Starting from
$\psi^{\pm}$ and $\phi_{\pm}$ given in (\ref{efunc1}) and
(\ref{efunc2}), we generate the following eigenfunctions
corresponding to the solution $Q$ given in (\ref{newsl}) through 
the iterated \BD,
\begin{eqnarray}
\Psi^{\pm} &=& \Gamma^{(I)} \Gamma \psi^{\pm} \, , \ \ \ \ \mbox{ at } 
\ \lambda = \lambda_0 = \sqrt{\eta^2- \frac{1}{4}\k_1^2} \, , 
\label{nwef1} \\[1ex]
\Phi_{\pm} &=& \Gamma^{(I)} \Gamma \phi_{\pm} \, ,  \ \ \ \
\mbox{ at } \lambda =0 \, , \label{nwef2}
\end{eqnarray}
where
\[
  \Gamma = \left[
    \begin{array}{cc}
\Lambda +a & b\\ \\ 
\overline{b} & \Lambda -a
    \end{array} \right] \, , \ \ \ \Gamma^{(I)} = \left[
    \begin{array}{cc}
\Lambda + a^{(I)} & b^{(I)} \\ \\ 
\overline{b^{(I)}} & \Lambda - \overline{a^{(I)}}
    \end{array} \right] \, ,
\]
where $\Lambda = \alpha \partial_y - \lambda$ for general $\lambda$.
\begin{lemma}
The eigenfunctions $\Psi^{\pm}$ and $\Phi_{\pm}$ defined in
(\ref{nwef1}) and (\ref{nwef2}) have the representations,
\begin{eqnarray}
    \Psi^{\pm} &=& \pm i \la_0 \k_1 \eta^{-1}\sqrt{c_0^+c_0^-} 
e^{ir_2t} [|v_1|^2+|v_2|^2]^{-1} \nonumber \\
& & \quad \quad \times \left[ \begin{array}{c}
-q_c \bigg [(\lambda_0 - a^{(I)}) \overline{v_2}+\eta \widetilde{b^{(I)}} 
\overline{v_1}\bigg ] \cr \cr
\eta \bigg [-\eta \overline{\widetilde{b^{(I)}}}\overline{v_2} - (\lambda_0 + 
\overline{a^{(I)}})\overline{v_1}\bigg ]\cr \end{array}\right] \, ,
\label{rnef1} \\
\Phi_{\pm} &=& \pm i \frac{1}{4}\al \k_2 \sqrt{c^0_+c^0_-} 
e^{i\xi_1^0x+ir_2t}[|W_1|^2+|W_2|^2]^{-1}
\left[ \begin{array}{c}
-q_c \widetilde{\Sg}_1 \cr \cr \eta \widetilde{\Sg}_2 \cr
\end{array} \right] \, , \label{rnef2}
\end{eqnarray}
where $b^{(I)}=-q_c \widetilde{b^{(I)}}$, and 
\begin{eqnarray*}
\widetilde{\Sg}_1 &=& 2\overline{W_1} (W_1^+W_1^-) +\overline{W_2^+} (W_1^+W_2^-) 
 +\overline{W_2^-} (W_1^-W_2^+)\ , \\
\widetilde{\Sg}_2 &=& 2\overline{W_2} (W_2^+W_2^-) +\overline{W_1^+} (W_2^+W_1^-) 
 +\overline{W_1^-} (W_2^-W_1^+)\ .
\end{eqnarray*}
If we take $r_2$ to be real [in the Melnikov vectors, $r_2$
appears in the form $r_2-r_1=2(\left|q_c \right|^2 - \omega^2)$], then
\begin{equation}
  \Psi^{\pm} \to 0 \, , \ \ \ \  \Phi_{\pm} \to 0 \, , \, 
  \hbox{ as } t \to \pm \infty \, .
\label{anef}
\end{equation}
\end{lemma}
Next we generate eigenfunctions solving the corresponding
congruent Lax pair (\ref{CLP1}, \ref{CLP2}) with the potential $Q$, through
the iterated \BD ~and the symmetry transformation (\ref{IT}) \cite{Li00}.
\begin{lemma}
Under the replacement
\begin{eqnarray}
& &  \alpha \longrightarrow - \alpha \, \quad (\mbox{then}\ \vth_2 
\longrightarrow \pi - \vth_2 ), \quad r_1  \longrightarrow -r_2, \nonumber \\
& & r_2  \longrightarrow -r_1, \quad
    \hvth \longrightarrow \hvth +\pi -2 \vth_2 \, , \quad 
    \hat{\rho} \longrightarrow - \hat{\rho}\, , \label{rptr}
\end{eqnarray}
the potentials are transformed as follows,
\begin{eqnarray*}
& & Q  \longrightarrow Q \, , \\
& & R_1 \longrightarrow -R_2 \, , \\
& & R_2 \longrightarrow -R_1 \, .
\end{eqnarray*}
\label{rppl}
\end{lemma}
The eigenfunctions $\Psi^{\pm }$ and $\Phi_{\pm}$ given in
(\ref{rnef1}) and (\ref{rnef2}) depend on the variables in the
replacement (\ref{rptr}):
\begin{eqnarray*}
  \Psi^{\pm} &=& \Psi^{\pm} (\alpha , r_1 , r_2 , \hvth , \hrho ) \, , \\
  \Phi_{\pm} &=& \Phi_{\pm} (\alpha , r_1 , r_2 , \hvth , \hrho ) \, .
\end{eqnarray*}
Under replacement (\ref{rptr}), $\Psi^{\pm}$ and $\Phi_{\pm}$ are 
transformed into
\begin{eqnarray}
  \widehat{\Psi}^{\pm} &=& \Psi^{\pm}
     (- \alpha , -r_2 , -r_1 , \hvth + \pi -2 \vth_2 , - \hrho )
     \, , \label{cref1}\\
  \widehat{\Phi}_{\pm} &=& \Phi_{\pm}
     (- \alpha , -r_2 , -r_1 , \hvth + \pi -2 \vth_2 , - \hrho )
     \, . \label{cref2}
\end{eqnarray}
\begin{lemma}
$\widehat{\Psi}^{\pm}$ and $ \widehat{\Phi}_{\pm}$ solve the
congruent Lax pair (\ref{CLP1}, \ref{CLP2}) at $(Q, \overline{Q}, 
R_1, R_2; \lambda_0)$ and $(Q, \overline{Q}, R_1, R_2; 0)$, respectively.
\label{ccor}
\end{lemma}
Notice that as a function of $\eta$, $\xi^0_1$ has two (plus and
minus) branches.  In order to construct Melnikov vectors, we need 
to study the effect of the replacement $\xi^0_1 \longrightarrow
-\xi^0_1$.
\begin{lemma}
Under the replacement
\begin{equation}
    \xi^0_1 \longrightarrow - \xi^0_1 \,  \ \ (\mbox{then}\ 
\vth_2 \longrightarrow -\vth_2 ), \ \  
   \hvth \longrightarrow \hvth + \pi -2 \vth_2 , \ \ 
\hrho \longrightarrow - \hrho  ,
\label{krptr}  
\end{equation}
the potentials are invariant,
\[
Q \longrightarrow Q \, , \quad R_1 \longrightarrow R_1 \, , \quad 
R_2 \longrightarrow R_2 \, .
\]
\label{lektr}
\end{lemma}
The eigenfunction $\Phi_{\pm}$ given in (\ref{rnef2}) depends on
the variables in the replacement (\ref{krptr}):
\begin{displaymath}
  \Phi_{\pm} = \Phi_{\pm} (\xi^0_1 , \hvth , \hrho ) \, .
\end{displaymath}
Under the replacement (\ref{krptr}), $\Phi_{\pm}$ is transformed
into
\begin{equation}
  \widetilde{\Phi}_{\pm} = \Phi_{\pm} 
  (-\xi^0_1 , \hvth + \pi -2 \vth_2, -\hrho )\ .
\label{kfr}
\end{equation}
\begin{lemma}
$\widetilde{\Phi}_{\pm}$ solves the Lax pair (\ref{LP1},\ref{LP2}) at
$(Q , \overline{Q} , R_1 , R_2 \, ; \, 0)$.
\label{kcor}
\end{lemma}
In the construction of the Melnikov vectors, we need to replace
$\Phi_{\pm}$ by $ \widetilde{\Phi}_{\pm}$ to guarantee the
periodicity in $x$ of period $\frac{2\pi}{\k_1}$.

The Melnikov vectors for the Davey-Stewartson~II equations are
given by,
\begin{eqnarray}
  \U^{+} &=& \left( \begin{array}{c}
\Psi^{+}_2 \widehat{\Psi}^{+}_2 \\[1ex]
\Psi^{+}_1 \widehat{\Psi}^{+}_1    
    \end{array}\right)^- +S \left(
\begin{array}{c}
     \Psi^{+}_2 \widehat{\Psi}^{+}_2 \\[1ex]
\Psi^{+}_1 \widehat{\Psi}^{+}_1    
    \end{array} \right) \, , \label{mv1}\\[2ex]
  \U_{+} &=& \left(
\begin{array}{c}
      \widetilde{\Phi}_{+}^{(2)} \widehat{\Phi}_{+}^{(2)} \\[1ex]
\widetilde{\Phi}_{+}^{(1)} \widehat{\Phi}_{+}^{(1)}    
\end{array}\right)^- +S \left(
\begin{array}{c}
     \widetilde{\Phi}_{+}^{(2)} \widehat{\Phi}_{+}^{(2)} \\[1ex]
\widetilde{\Phi}_{+}^{(1)} \widehat{\Phi}_{+}^{(1)}    
    \end{array} \right) \, , \label{mv2}
\end{eqnarray}
where ``--'' denotes complex conjugate, and $S = \left ( 
\begin{array}{lr} 0 & 1 \\ 1 & 0 \end{array} \right )$. In fact, the even 
parts of $\U^{+}$ and $\U_{+}$ are the Melnikov vectors in our phase space.
Nevertheless, the Melnikov integral formulas end up the same, as shown in 
\cite{Li00}. For simplicity, we just use $\U^{+}$ and $\U_{+}$.

\subsection{Melnikov Analysis}

The main difficulty in a rigorous Melnikov measurement is due to the lack 
of global well-posedness. The main idea in resolving this difficulty is to 
iterate the small time interval in local well-posedness by virtue of the fact 
that the unperturbed homoclinic orbit is a classical solution.

Let $p$ be any point on $W^u_\Pi(Q_\e)$, the unstable curve of $Q_\e$ in $\Pi$.
By the Unstable Fiber Theorem \ref{UFT}, $\F^u_p$ is $C^1$ in $\e$ for 
$\e \in [0,\e_0)$, $\e_0 > 0$; thus, there are two points $q_\e(0)$ and 
$q_0(0)$ on the unstable fibers $\F^u_p$ and $\F^u_p|_{\e =0}$, such that
\[
\| q_\e(0)-q_0(0)\|_n \leq C_n^{(1)} \e \ , \quad (n \geq 2)\ .
\]
The key point here is that $\mathcal{F}^u_p\subset H^s$ for any fixed 
$s\geq 2$. The expression of the unperturbed homoclinic orbit $q_0(t)$ 
has been given in (\ref{newsl}) which represents a classical solution 
to the DSII. Let
\[
D^*_s = \sup_{t \in (-\infty, +\infty)} \{ \ \| q_0(t) \|_s \ \} \ , \quad 
(s \geq 2)\ .
\]
By the Local Well-Posedness Theorem \ref{locwp}, there exists 
$\tau = \tau(D^*_n)>0$, such that
\[
\| q_\e(t)-q_0(t)\|_n \leq C_n^{(2)} \e \ , \quad t \in [0,\tau ]\ ,
\]
where $C_n^{(2)} = C_n^{(2)} (D^*_{n+4})$. There is an integer $N>0$ such 
that 
\[
q_0(N\tau ) \in W^{cs}_n|_{\e =0}\ ,
\]
where $W^{cs}_n$ is given by the Center-Stable Manifold Theorem \ref{CSM}. 
Iterating the Local Well-Posedness Theorem $N$ times, one gets
\[
\| q_\e(t)-q_0(t)\|_n \leq C_n^{(3)} \e \ , \quad t \in [0,N\tau ]\ ,
\]
where $C_n^{(3)} = C_n^{(3)} (D^*_{n+4})$. Our goal is to determine when 
$q_\e(N\tau ) \in W^{cs}_n$ through Melnikov measurement. The two Melnikov 
vectors $\U^+$ and $\U_+$ (\ref{mv1})-(\ref{mv2}) are transversal to 
$W^{cs}_n$. There is a unique point $\hq_\e(N\tau) \in W^{cs}_n$ such that
\[
q_\e(N\tau )- \hq_\e(N\tau) \in \ \mbox{span}\ \{ \U^+,\U_+ \}\ ;
\]
thus, $\hq_\e(N\tau) \in W^{cs}_{n+4}$. By the Center-Stable Manifold 
Theorem \ref{CSM},
\[
\| \hq_\e(N\tau) -q_0(N\tau )\|_n \leq C_n^{(4)} \e \ , 
\]
where $C_n^{(4)} = C_n^{(4)} (D^*_{n+4})$. Thus
\[
\| q_\e(N\tau )-\hq_\e(N\tau)\|_n \leq C_n \e \ ,
\]
where $C_n = C_n (D^*_{n+4})$. To determine when $q_\e(N\tau )=\hq_\e(N\tau)$, 
one can define the signed distances
\[
d_1 = \lag \U^+, \vec{q_\e}(N\tau ) - \vec{\hq_\e}(N\tau ) \rag\ , \quad 
d_2 = \lag \U_+, \vec{q_\e}(N\tau ) - \vec{\hq_\e}(N\tau ) \rag\ , 
\]
where $\vec{q}=(q, \bq)^T$, and 
\[
\lag A, B \rag = \int_0^{2\pi/\k_2} \int_0^{2\pi/\k_1} \{ \overline{A_1}
B_1 + \overline{A_2}B_2 \} \ dx dy \ .
\]
The rest of the derivation for Melnikov integrals is completely standard.
For details, see e.g. \cite{LM97} \cite{LMSW96}.
\[
d_k = \e M_k +o(\e)\ , \quad k=1,2,
\]
where
\[
M_1 = \int_{-\infty}^\infty \lag \U^+, G \rag \ dt \ , \quad 
M_2 = \int_{-\infty}^\infty \lag \U_+, G \rag \ dt \ , 
\]
where $G=(f,\bar{f})^T$, $f=\Dl Q -\al Q +\be$. That is,
\begin{eqnarray*}
M_1 &=& \int_{-\infty}^\infty \int_0^{2\pi/\k_2} \int_0^{2\pi/\k_1} 
\ \mbox{Re}\ \{ (\Psi^{+}_2 \widehat{\Psi}^{+}_2)f + 
(\Psi^{+}_1 \widehat{\Psi}^{+}_1)\bar{f} \} \ dxdydt\ , \\
M_2 &=& \int_{-\infty}^\infty \int_0^{2\pi/\k_2} \int_0^{2\pi/\k_1} 
\ \mbox{Re}\ \{ (\widetilde{\Phi}_{+}^{(2)} \widehat{\Phi}_{+}^{(2)})f + 
(\widetilde{\Phi}_{+}^{(1)} \widehat{\Phi}_{+}^{(1)})\bar{f} \} \ dxdydt\ , 
\end{eqnarray*}
where $\eta =\om$, and we divide $\Psi^{+}$ by the constant $i\la_0 \k_1 
\sqrt{c_0^+c_0^-} e^{i\ga/2}$, and $\Phi_+$ by $\frac{1}{4} i \al \k_2 \eta 
\sqrt{c^0_+c^0_-} e^{i\ga/2}$. It has been verified numerically that 
multiplication of $\Psi^{+}$ and $\Phi_+$ by a complex constant leads to 
equivalent results. It turns out that
\[
M_j = M_j^{(1)}+\al M_j^{(2)}+\be \cos \ga M_j^{(3)}+ \be \sin \ga 
M_j^{(4)} \ , \quad (j=1,2)\ ,
\]
where $M_j^{(l)}=M_j^{(l)}(\om, \Dl \rho )$, ($j=1,2; 1 \leq l 
\leq 4$), $\Dl \rho = \hat{\rho} + i\al \k_2 \xi_1^0 \k_1^{-1}\la_0^{-1}
\rho$, $\htau = i\al \k_2 \xi_1^0 \k_1^{-1}\la_0^{-1}\tau +\Dl \rho $.

$M_j = 0$ ($j=1,2$) imply that 
\begin{eqnarray}
\al = \al (\om,\Dl \rho, \ga) &=& \bigg \{ M_1^{(1)} [\cos \ga 
M_2^{(3)}+\sin \ga M_2^{(4)}] \nonumber \\
& & - M_2^{(1)} [\cos \ga 
M_1^{(3)}+\sin \ga M_1^{(4)}] \bigg \} \nonumber \\
& & \times \bigg \{ M_2^{(2)} [\cos \ga 
M_1^{(3)}+\sin \ga M_1^{(4)}] \nonumber \\
& & - M_1^{(2)} [\cos \ga 
M_2^{(3)}+\sin \ga M_2^{(4)}] \bigg \}^{-1}\ ,  \label{alv}\\
\be = \be (\om,\Dl \rho, \ga) &=& [ M_1^{(1)}M_2^{(2)}-M_2^{(1)}
M_1^{(2)}] \nonumber \\
& & \times \bigg \{ M_1^{(2)} [\cos \ga 
M_2^{(3)}+\sin \ga M_2^{(4)}] \nonumber \\
& & - M_2^{(2)} [\cos \ga M_1^{(3)}+\sin \ga M_1^{(4)}] \bigg \}^{-1}\ . 
\label{bev} 
\end{eqnarray}
\begin{theorem}
There exists $\epsilon_0>0$, such that 
for any $\e \in (0, \e_0)$, there exists a domain $\mathcal{D}_\e 
\subset \mathbb{R}^+ \times \mathbb{R}^+ \times \mathbb{R}^+$ where 
$\om$ satisfies the constraint (\ref{cstr1}) or (\ref{cstr2}), 
and $\al \om < \be$. For any $(\alpha ,\beta, \omega) \in \mathcal{D}_\e$,
there exists another orbit in $W^u(Q_\epsilon )\cap W^{cs}_n$ other 
than the unstable curve $W^u_\Pi (Q_\epsilon )$ of $Q_\e$ in $\Pi$, for 
the perturbed DSII (\ref{PDS}).
\end{theorem}

Proof. The zeros of $M_j$ ($j=1,2$) are given by (\ref{alv}) and 
(\ref{bev}). We need $\al >0$ and $\be >0$ which define a region in 
the external parameter space, parametrized by $\Dl \rho$ and $\ga$. 
Then the theorem follows from the implicit function theorem. Q.E.D.

For example, when $\k_1 =1$ and $\k_2 =\sqrt{2}$, 
\[
\al (\frac{\sqrt{2}}{2}+0.11, 1.1, \frac{\pi }{2}) = 5.645\ , 
\quad \be (\frac{\sqrt{2}}{2}+0.11, 1.1, \frac{\pi }{2}) = 11.336 \ .
\]

%% file: appendix.tex
\section{Appendix}

The main obstacle toward proving the existence of a homoclinic orbit for 
the perturbed DSII (\ref{PDS}) comes from a technical difficulty in the 
normal form transform \cite{Li01b}. In this appendix, we will present 
the difficulty.

\subsection{The Technical Difficulty in the Normal Form Transform}

To locate a homoclinic orbit to $Q_\epsilon$ \eqref{Qec}, we need 
to estimate the size of the local stable manifold of $Q_\epsilon$.
The size of the variable $J$ is of order $\mathcal{O}(\sqrt{\epsilon})$. 
The size of the variable $\theta$ is of order $\mathcal{O}(1)$.
To be able to track a homoclinic orbit, we need the size of the 
variable $f$ to be of order $\mathcal{O}(\epsilon^\mu)$, $\mu <1$. 
Such an estimate can be achieved, if the quadratic term 
$\mathcal{N}_2$ in \eqref{cc3} can be removed through a 
normal form transformation. In fact, it is enough to remove its leading 
order part
\[
\tilde{\N}_2 = 2\om \bigg [ \Delta^{-1}\sq |f|^2
+f\Delta^{-1} \sq (f+\bar f) 
- \lag f\Delta^{-1}\sq(f+\bar f)\rag \bigg ].
\]
That is, our goal is to find a normal form transform $g = f + K(f,f)$ 
where $K$ is a bilinear form, that transforms the equation
\[
f_t=L_\epsilon f-i\tilde{\mathcal{N}}_2,
\]
into an equation with a cubic nonlinearity
\[
g_t=L_\epsilon g+\mathcal{O}(\| g\|^3_s),\quad 
(s\geq 2),
\]
where $L_\epsilon$ is given in (\ref{cc3}). In terms of Fourier transforms,
\[
f=\sum_{k\neq 0}\hat{f}(k)e^{ik\cdot \xi },\quad 
\bar{f}=\sum_{k\neq 0}\overline{\hat{f}(-k)}e^{ik\cdot \xi }\ ,
\]
where $k=(k_1,k_2) \in \mathbb{Z}^2$, $\xi = (\k_1x, \k_2y)$. The terms in 
$\tilde{\mathcal{N}}_2$ can be written as
\[
\Dl^{-1}\sq |f|^2 = \frac {1}{2} \sum_{k+\ell \neq 0} a(k+\ell) 
\bigg [ \hat{f}(k)\overline{\hat{f}(-\ell )}+\hat{f}(\ell )
\overline{\hat{f}(-k)} \bigg ] e^{i(k+ \ell ) \cdot \xi }\ ,
\]
\[
f\Delta^{-1} \sq f - \lag f\Delta^{-1}\sq f \rag  
= \frac {1}{2} \sum_{k+\ell \neq 0} [a(k)+a( \ell )] \hat{f}(k)\hat{f}(\ell)
e^{i(k+ \ell ) \cdot \xi }\ ,
\]
\[
f\Delta^{-1} \sq \bar{f} - \lag f\Delta^{-1}\sq \bar{f} \rag  
= \frac {1}{2} \sum_{k+\ell \neq 0} \bigg [ a(\ell)\hat{f}(k)
\overline{\hat{f}(-\ell )}+a(k)\hat{f}(\ell)
\overline{\hat{f}(-k)}\bigg ] e^{i(k+ \ell ) \cdot \xi }\ ,
\]
where
\[
a(k) = \frac {k_1^2 \k_1^2 - k_2^2 \k_2^2}{k_1^2 \k_1^2 + k_2^2 \k_2^2}\ .
\]
We will search for a normal form transform of the 
general form,
\[
g=f+K(f,f),
\]
where
\begin{eqnarray*}
K(f,f) &=& \sum_{k+\ell\neq 0}\left[
\hat{K}_1(k,\ell)\hat{f}(k) 
\hat{f}(\ell)+\hat{K}_2(k,\ell)\hat{f}(k)\overline{\hat{f}(-\ell)}\right.\\
& &\quad \left.+\hat{K}_2(\ell,k)\overline{\hat{f}(-k)}\hat{f}(\ell)+
\hat{K}_3(k,\ell)\overline{\hat{f}(-k)}\overline{\hat{f}(-\ell)}\right]
e^{i(k+\ell)x},
\end{eqnarray*}
where $\hat{K}_j(k,\ell)$, $(j=1,2,3)$ are the unknown coefficients to 
be determined, and $\hat{K}_j(k,\ell)=\hat{K}_j(\ell,k)$, $(j=1,3)$.
To eliminate the quadratic terms, we first need to set
\[
iL_\epsilon K(f,f)-iK(L_\epsilon 
f,f)-iK(f,L_\epsilon f)=\tilde{\mathcal{N}}_2,
\]
which takes the explicit form:
\begin{eqnarray}
& &(\sigma_1+i\sigma)\hat{K}_1(k,\ell)+B(\ell)\hat{K}_2(k,\ell)+B(k)
\hat{K}_2(\ell,k)\nonumber \\
& & \quad +B(k+\ell)\overline{\hat{K}_3(k,\ell)}=\frac{1}{2\om} 
[B(k)+B(\ell)],\label{nore1}\\
& &-B(\ell)\hat{K}_1(k,\ell)+(\sigma_2+i\sigma)\hat{K}_2(k,\ell)+ B(k+\ell)
\overline{\hat{K}_2(\ell,k)}\nonumber \\
& & \quad +B(k)\hat{K}_3(k,\ell)=\frac{1}{2\om} 
[B(k+\ell)+B(\ell)],\label{nore2}\\
& &-B(k)\hat{K}_1(k,\ell)+B(k+\ell)
\overline{\hat{K}_2(k,\ell)}+(\sigma_3+i\sigma)\hat{K}_2(\ell,k)\nonumber \\
& & \quad +B(\ell)\hat{K}_3(k,\ell)= \frac{1}{2\om} 
[B(k+\ell)+B(k)],\label{nore3}\\
& &B(k+\ell) \overline{\hat{K}_1(k,\ell)}-B(k)\hat{K}_2(k,\ell) 
-B(\ell)\hat{K}_2(\ell,k)\nonumber \\
& & \quad +(\sigma _4+i\sigma )\hat{K}_3(k,\ell)=0,\label{nore4}
\end{eqnarray}
where $B(k)=2\om^2 a(k)$, and 
\begin{eqnarray*} 
& & \sg = \e \bigg [ \al -2(k_1\ell_1\k_1^2+k_2\ell_2\k_2^2)\bigg ]\ , \\
& & \sg_1= 2(k_2\ell_2\k_2^2-k_1\ell_1\k_1^2) +B(k+\ell)-B(k)-B(\ell)\ , \\
& & \sg_2= 2[(k_2+\ell_2)\ell_2\k_2^2-(k_1+\ell_1)\ell_1\k_1^2] 
+B(k+\ell)-B(k)+B(\ell)\ , \\
& & \sg_3= 2[(k_2+\ell_2)k_2\k_2^2-(k_1+\ell_1)k_1\k_1^2] 
+B(k+\ell)+B(k)-B(\ell)\ , \\
& & \sg_4= 2[(k_2^2+k_2\ell_2+\ell_2^2)\k_2^2-(k_1^2+k_1\ell_1+\ell_1^2)
\k_1^2] +B(k+\ell)+B(k)+B(\ell)\ . 
\end{eqnarray*}
Since these coefficients are even in $(k,\ell)$, we will search for 
even solutions, i.e.
\[
\hat{K}_j(-k,-\ell)=\hat{K}_j(k,\ell),\quad j=1,2,3.
\]

The technical difficulty in the normal form transform comes from not 
being able to answer the following two questions in solving the linear 
system (\ref{nore1})-(\ref{nore4}):
\begin{enumerate}
\item Is it true that for all $k, \ell \in \mathbb{Z}^2/\{ 0\}$, there 
is a solution ?
\item What is the asymptotic behavior of the solution as $k$ and/or 
$\ell \ra \infty$ ? In particular, is the asymptotic behavior like 
$k^{-m}$ and/or $\ell^{-m}$ ($m\geq 0$) ?
\end{enumerate}

\subsection{A Formal Calculation}

Formally conducting the calculation for the second measurement to locate 
a homoclinic orbit \cite{Li01b}, one gets the formulas
\[
M_j =0 \quad (j=1,2)\ , \quad \be \cos \ga = -\frac{ \al \om \Dl \ga}
{2\sin \frac{\Dl \ga}{2}}\ ,
\]
where $\Dl \ga =-4 (\vth_1 -\vth_2)$. Thus we have $\al = 1/\chi$,
\begin{eqnarray*}
\chi =\chi (\om, \Dl \rho) &=& (M_2^{(1)}M_1^{(4)}-M_1^{(1)}M_2^{(4)})^{-1}
\bigg[ M_1^{(2)}M_2^{(4)}-M_2^{(2)}M_1^{(4)} \\
& & -\om \Dl \ga [2\sin \frac{\Dl \ga}{2}]^{-1} (M_1^{(3)}M_2^{(4)}-
M_2^{(3)}M_1^{(4)})\bigg ]\ , \\
\be = \be (\om, \Dl \rho) &=& \bigg[ (\al \om \Dl \ga)^2
[2\sin \frac{\Dl \ga}{2}]^{-2}\\
& & +(M_1^{(4)})^{-2}[M_1^{(1)} + \al (
M_1^{(2)}-M_1^{(3)}\om \Dl \ga ( 2\sin \frac{\Dl \ga}{2})^{-1})]^2 
\bigg ]^{1/2}\ .
\end{eqnarray*}

For example, when $\k_1 =1$ and $\k_2 =\sqrt{2}$, 
\[
\chi (\frac{\sqrt{2}}{2}+0.11, 1.1)=0.4326\ .
\]